\documentclass[11pt]{article}

\usepackage{enumerate}
\usepackage{amsmath,xspace,amssymb,mathrsfs}

\input xy
\xyoption{all}
\xyoption{2cell}
\UseAllTwocells
\CompileMatrices

\title{Icons}
\author{Stephen Lack%
\thanks{The support of the Australian Research Council and
DETYA is gratefully acknowledged.}
\\ School of Computing and Mathematics\\
University of Western Sydney\\
Locked Bag 1797 Penrith South DC NSW 1797\\
Australia
\\ email: {\tt s.lack@uws.edu.au}
}
\date{}

\renewcommand{\phi}{\varphi}

\newcommand{\A}{{\ensuremath{\mathscr A}}\xspace}
\newcommand{\B}{{\ensuremath{\mathscr B}}\xspace}
\newcommand{\C}{{\ensuremath{\mathscr C}}\xspace}
\newcommand{\D}{{\ensuremath{\mathscr D}}\xspace}
\newcommand{\G}{{\ensuremath{\mathscr G}}\xspace}
\newcommand{\HH}{{\ensuremath{\mathscr H}}\xspace}
\newcommand{\K}{{\ensuremath{\mathscr K}}\xspace}
\newcommand{\T}{{\ensuremath{\mathscr T}}\xspace}
\newcommand{\V}{{\ensuremath{\mathscr V}}\xspace}
\newcommand{\W}{{\ensuremath{\mathscr W}}\xspace}
\newcommand{\N}{{\ensuremath{\mathbb N}}\xspace}

\newcommand{\sS}{{\ensuremath{\mathscr S}}\xspace}
\newcommand{\cC}{{\ensuremath{\mathscr C}}\xspace}

\newcommand{\cat}{\ensuremath{\textnormal{\bf Cat}_\textnormal{\bf 1}}\xspace}
\newcommand{\Cattwo}{\ensuremath{\textnormal{\bf Cat}_\textnormal{\bf 2}}\xspace}
\newcommand{\bicat}{\ensuremath{\textnormal{\bf Bicat}_\textnormal{\bf 2}}\xspace}
\newcommand{\homm}{\ensuremath{\textnormal{\bf Hom}_\textnormal{\bf 2}}\xspace}
\newcommand{\nhom}{\ensuremath{\textnormal{\bf nHom}_\textnormal{\bf 2}}\xspace}
\newcommand{\shom}{\ensuremath{\textnormal{\bf sHom}_\textnormal{\bf 2}}\xspace}
\newcommand{\Oplax}{\ensuremath{\textnormal{\bf Oplax}_\textnormal{\bf 2}}\xspace}

\newcommand{\LaxDbl}{{\textnormal{\bf LaxDbl}}\xspace}
\newcommand{\Cat}{{\textnormal{\bf Cat}}\xspace}
\newcommand{\Hom}{\textnormal{\bf Hom}\xspace}

\newcommand{\CG}{\ensuremath{\textnormal{\bf Cat-Gph}_\textnormal{\bf 2}}\xspace}
\newcommand{\talg}{\textnormal{\bf $T$-Alg}\xspace}
\newcommand{\talgl}{\ensuremath{\textnormal{\bf $T$-Alg}_{\textnormal{l}}}\xspace}

\newcommand{\Set}{{\textnormal{\bf Set}}\xspace}
\newcommand{\Grp}{{\textnormal{\bf Grp}}\xspace}
\newcommand{\Mod}{\textnormal{\bf Mod}\xspace}
\newcommand{\MonCat}{\textnormal{\bf MonCat}\xspace}
\newcommand{\Lex}{\textnormal{\bf Lex}\xspace}
\newcommand{\lfp}{\textnormal{\bf LFP}\xspace}
\newcommand{\n}{\mathbf{n}}
\newcommand{\op}{\ensuremath{{}^{\textnormal{op}}}\xspace}

\newcommand{\two}{\ensuremath{\mathbf{2}}\xspace}

\renewcommand{\t}{\times}
\newcommand{\ot}{\otimes}

\newtheorem{theorem}{Theorem}[section]

\newtheorem{proposition}[theorem]{Proposition}
\newtheorem{preremark}[theorem]{Remark}
\newenvironment{remark}{\begin{preremark}\rm}{\end{preremark}}
\newtheorem{predefinition}[theorem]{Definition}
\newenvironment{definition}{\begin{predefinition}\rm}{\end{predefinition}}

\newcommand{\proof}{\noindent{\sc Proof:}\xspace}
\def\endproof{~\hfill$\Box$\vskip 10pt}

\begin{document}

\label{firstpage}
\maketitle

\begin{abstract}
Categorical orthodoxy has it that collections of ordinary mathematical 
structures such as groups, rings, or spaces, form {\em categories} (such as the 
category of groups); collections of 1-dimensional categorical structures, such
as categories, monoidal categories, or categories with finite 
limits, form {\em 2-categories}; and collections of 2-dimensional categorical 
structures, such as 2-categories or bicategories, form {\em 3-categories}.

We describe a useful way in which to regard bicategories as objects of a 
2-category. This is a bit surprising both for technical and for conceptual
reasons. The 2-cells of this 2-category are the crucial new ingredient;
they are the icons of the title. These can be thought of as ``the
oplax natural transformations whose components are identities'', but 
we shall also give a more elementary description.

We describe some properties of these icons,
and give applications to monoidal categories, to  2-nerves of bicategories, 
to 2-dimensional Lawvere theories, and to bundles of bicategories.
\end{abstract}

\section{Introduction}

For any particular mathematical structure, there is a category whose
objects are instances of that structure, and whose morphisms are the 
structure-preserving maps. For the structure of group, we have the category
\Grp of groups and group homomorphisms. Similarly, for the structure of 
category itself, we have the category \cat of categories and functors. But 
in this case, as well as functors, there
are also natural transformations between functors, and the category
of categories naturally becomes a 2-category called \Cattwo or just \Cat. 

What happens when the mathematical structure in question is a 2-dimensional
one, such as a bicategory? Once again, one could use the 
structure preserving maps to obtain a category of bicategories (although
there is a certain amount of choice in what exactly is meant by 
``structure-preserving''), but now there 
are two further levels of structure, which altogether give a 
3-dimensional structure called a tricategory \cite{GPS} or weak 3-category.
In particular, there is a tricategory called \Hom \cite{FibBic,GPS} whose
objects are the bicategories. 
There is no doubt that this is a fine and excellent structure.
But my goal here is to argue for the merits of a different context in which
to consider bicategories. This different context will be just 2-dimensional,
so that we get a 2-category of bicategories. This is not done by 
merely throwing away the 3-cells of \Hom --- for one thing this would not
give a 2-category or even a bicategory, but even if it did it would conceptually
not be the right thing to do. Rather, we introduce a new type of 2-cell,
which we call an {\em icon}, and offer the following propaganda for 
these new 2-cells.

\paragraph{They provide a relatively simple context in which to study 
bicategories.}

The 3-dimensional structure \Hom is, as mentioned above, a 
tricategory \cite{GPS}, and
tricategories are rather complex structures. Although for some purposes
one really needs to use this full tricategory, if one can get by with 
just a 2-category of bicategories life is very much simpler. For example, 
this allows bicategories to be treated using the techniques of 2-dimensional 
universal algebra \cite{BKP}, which is very much more developed than 
3-dimensional universal algebra \cite{Power:3BKP}.

\paragraph{They arise naturally  if one thinks of 
2-categories or bicategories as ``many-object'' monoidal categories.}
One-object categories are essentially equivalent to monoids. If you think 
of monoids as
one-object categories, then a homomorphism of monoids is just a functor.
Thus the category of monoids is a full subcategory of the category \cat
of categories and functors; this does not use the 2-category structure
\Cat built upon \cat.

Similarly, one-object bicategories are essentially equivalent to  
monoidal categories, and if you think of monoidal categories as one-object 
bicategories, then structure preserving morphisms of bicategories 
are the same as structure preserving morphisms of monoidal categories.
But between the structure-preserving morphisms of monoidal categories
there are things called monoidal natural transformations, and these 
are precisely the icons (in the case of one-object bicategories).

\paragraph{They facilitate the construction of 2-nerves of bicategories.}
The nerve construction is of great importance in algebraic topology and 
elsewhere. It associates to every category a simplicial set; in fact it
provides a full embedding of the category \cat of categories and functors
into the category $[\Delta\op,\Set]$ of simplicial sets. This construction
is a formal consequence of the fact that there is an inclusion of $\Delta$ 
in \cat. The 2-nerve 
construction associates to every bicategory a simplicial category. We can
obtain this as a formal consequence of the fact that $\Delta$
embeds in a suitably defined 2-category of bicategories, with icons
as 2-cells, as was proved in  \cite{2-nerve}.

\paragraph{They arise in connection with 2-dimensional Lawvere theories.}
Lawvere theories are used to describe one-dimensional algebraic structures,
such as groups, rings, or Lie algebras. 2-dimensional Lawvere theories 
\cite{Power:enriched-lawvere,2dLawvere} are
used to describe two-dimensional algebraic structures, such as monoidal
categories or categories with finite limits. Formally, a Lawvere theory is
a particular type of category, and a 2-dimensional Lawvere theory is a 
special type of 2-category. A morphism of 2-dimensional Lawvere theories
is a special type of 2-functor, and a transformation between such morphisms
is a special type of icon. These transformations play 
a key role in a notion called {\em flexibility}, which distinguishes structures
such as {\em monoidal categories} which transport along equivalences, 
to structures such as {\em strict monoidal categories}, which do not.

\paragraph{They allow the use of notions internal to a 2-category,
such as equivalence and fibration.}
Various notions such as equivalence and fibration can be defined internally
to any 2-category; in particular this can be done in the various 
2-categories of bicategories defined using icons. Equivalence in these 
2-categories is quite a strict notion, but it can be proved that any
bicategory is equivalent, in a suitable 2-category of bicategories with 
icons as 2-cells, to
a strict bicategory (that is, a 2-category) \cite{2-nerve}. 

On the other hand 
a fibration in the various 2-categories of bicategories can be thought of as a 
sort of ``bundle of bicategories'' \cite{HL}.

\section{Bicategorical preliminaries}
\label{sect:background}

As indicated above, and as first observed in \cite{bicategories}, 
there are various possible notions of morphism of 
bicategory. First one could consider those which literally preserve
all of the structure, so that for instance, given composable maps $f$
and $g$ we have $F(gf)=F(g)F(f)$. These are called {\em strict homomorphisms}.
In real life, these strict homomorphisms are rare; much more common are
{\em homomorphisms}, which preserve composition and identities only 
up to isomorphisms, such as $\phi_{g,f}:F(g)F(f)\cong F(gf)$ and 
$\phi^0_A:1_{FA}\cong F1_A$, which are themselves required to satisfy 
certain coherence conditions. When the $\phi^0_A$ are identities, so that the
identities at least are preserved strictly, a homomorphism is said to 
be {\em normal}. There is also a still more general notion, called a
{\em lax functor}, which involves comparisons $\phi_{g,f}:F(g)F(f)\to F(gf)$
and $\phi^0_A:1_{FA}\to F1_A$ which still satisfy the coherence conditions,
but are no longer required to be invertible.

There is a category of bicategories and lax functors, and this has various 
subcategories in which the morphisms are restricted to the homomorphisms,
the normal homomorphisms, or the strict homomorphisms, as the case may be.
Just as in the one-dimensional case, one could stop here, and consider
the category of bicategories and homomorphisms, or one of its variants,
and for some purposes this suffices. But just as one gains a deeper 
understanding of categories by regarding them as objects of a 2-category,
so one gains a deeper understanding of bicategories by regarding them
as the objects of a 3-dimensional structure called a tricategory \cite{GPS}.

Just as functoriality can be weakened, so can naturality, giving rise
to structures such as oplax natural transformations and pseudonatural
transformations, whose precise definitions are given in the following
section. These transformations serve as 2-cells, but there are also 
3-cells between them, called modifications. The most important 
3-dimensional category of bicategories is the tricategory \Hom
\cite{FibBic,GPS}, mentioned
in the introduction, which consists of the following levels of structure:
\begin{itemize}
\item bicategories as 0-cells 
\item homomorphisms between bicategories as 1-cells
\item pseudonatural transformations between homomorphisms as 2-cells
\item modifications between pseudonatural transformations as 3-cells.
\end{itemize}
Now the 2-cells of primary interest in this paper, called icons, are
not in general pseudonatural transformations. They are obtained by 
first generalizing from pseudonatural transformations to oplax natural
transformations, but then specializing to those whose components are 
identities, and so might perhaps be called {\em I}dentity {\em C}omponent 
{\em O}plax {\em N}atural transformations, whence the name {\em icon}.

As well as the various reasons for considering icons which were given 
in the introduction, we can now add:

\paragraph{They are the ``costrict'' transformations.} Among all the 
oplax natural transformations, one can identify the 2-natural transformations
abstractly via a notion of strictness. 
The dual notion of ``costrictness'' identifies exactly the icons.

Formally, the notion of strictness
is defined in a 2-dimensional structure called a {\em sesquicategory}
\cite{Street:categorical-structures}.
A sesquicategory \K consists of an underlying category $\K_1$, with 
a functor $\K:\K_1\op\t\K_1\to\Cat_1$ whose composite with the set-of-objects
functor $\Cat_1\to\Set$ is just the hom-functor $\K_1\op\t\K_1\to\Set$.
(There is also an alternative definition: sesquicategories are categories
enriched in \cat with the ``funny tensor product'': see 
\cite{Street:categorical-structures}.)

Spelling out the definition a little,  there are objects $A,B,C,\ldots$; 
1-cells like $f:A\to B$, and
2-cells $\alpha:f\to g:A\to B$. There is a strictly associative and unital
composition of 1-cells, and a strictly associative and unital {\em vertical}
composition of 2-cells (allowing the composite of $\alpha:f\to g$ and 
$\beta:g\to h$), and 2-cells can be composed on either side with 1-cells,
and this is also associative and unital.
On the other hand there is no specified {\em horizontal} composite of 
2-cells, as in 
$$\xymatrix{
\A \rtwocell^{F}_{G}{\alpha} & \B \rtwocell^{H}_{K}{\beta} & \C. }$$
Now we {\em can} form $\beta G$ and $H\alpha$, and then compose these, 
or alternatively we can form $K\alpha$ and $\beta F$ and compose them, 
and in each case we obtain a 2-cell from $HF$to $KG$, but there is no
need for the two resulting 2-cells to be equal, as in
$$\xymatrix @R1pc  {
\A \rtwocell^{F}_{G}{\alpha} & \B \ar@/^1pc/[r]^{H} & \C &&
\A \ar@/^1pc/[r]^{F} & \B \rtwocell^{H}_{K}{\beta} & \C \\
\A \ar@/_1pc/[r]_{G}  & \B \rtwocell^{H}_{K}{\beta} & \C &=&
\A \rtwocell^{F}_{G}{\alpha} & \B \ar@/_1pc/[r]_{K} & \C. }$$
This equation is sometimes called the ``middle-four interchange'' law.
A 2-category is precisely a sesquicategory in which this middle-four
interchange law holds for all $\alpha$ and $\beta$.

\begin{definition}
A 2-cell $\beta$ in a sesquicategory is {\em strict} if the equation above
holds for all $\alpha$.  
Dually, $\alpha$ is costrict if the equation holds for all $\beta$.
\end{definition}

\section{The definition}
\label{sect:definition}

Bicategories and lax functors form a category: in particular, 
composition of lax functors is strictly associative. There is also a subcategory
consisting only of the homomorphisms. Our goal is to make both 
of these into 2-categories by introducing 2-cells.

The associativity and identity isomorphisms $h(gf)\cong h(gf)$ and
$f1\cong f\cong 1f$ in a bicategory will rarely be mentioned explicitly, 
and will not be given a particular name. A lax functor $F:\A\to\B$ involves 
{\em lax functoriality constraints} $\phi_{g,f}:Fg.Ff\to F(gf)$ and 
{\em identity constraints} $\phi^0_A:1_{FA}\to F1_A$, usually abbreviated to 
$\phi$ and $\phi^0$. In the case of a lax functor called $G$ we use $\psi$
and $\psi^0$.

Consider bicategories \A and \B, and lax functors $F,G:\A\to\B$. 
The data for an {\em oplax natural transformation} from $F$ to $G$ consists of 
a 1-cell $\alpha A:FA\to GA$ in \B for every object $A\in\A$, and a 2-cell
$$\xymatrix @R1pc {
FA \ar[r]^{Ff} \ar[dd]_{\alpha A} & FB \ar[dd]^{\alpha B} \\
{}\rtwocell\omit{~~\alpha f} & {} \\
GA \ar[r]_{Gf} & GB }$$
in \B for every 1-cell $f:A\to B$ in \A. These are subject to three conditions.

\begin{enumerate}[(ON1)]
\addtocounter{enumi}{-1}
\item 
$\xymatrix@R1pc@C3pc{%
      FA \ruppertwocell<5>^{Ff}{<3>~~\alpha f} \ar[dd]_{\alpha A} &
      FB \ar[dd]^{\alpha B} &&
      FA \ar[dd]_{\alpha A} \rtwocell^{Ff}_{Fg}{~~F\rho} & FB \ar[dd]^{\alpha B} \\
      &&=&& \\
GA \rtwocell^{Gf}_{Gg}{~~G\rho} & GB && 
GA \rlowertwocell<-5>_{Gg}{<-3>~~\alpha g} & GB }$
\newline\noindent for all 2-cells $\rho:f\to g$ in \A
\item 
$\xymatrix @R1pc {%
FA \ar[r]^{Ff} \ar[dd]_{\alpha A} & FB \ar[r]^{Fg} \ar[dd]|(0.6){\alpha B} &
FC \ar[dd]^{\alpha C} && 
FA \ar[r]^{Ff} \ar[dd]_{\alpha A} \rrlowertwocell_{F(gf)}{~~\phi_{f,g}} & 
FB \ar[r]^{Fg} & FC \ar[dd]^{\alpha C} \\
{}\rtwocell\omit{~\alpha f} & {}\rtwocell\omit{~\alpha g} & {}&=&
\\
GA \ar[r]^{Gf} \rrlowertwocell_{G(gf)}{~~~\psi_{f,g}} & GB \ar[r]^{Gg} & GC &&
GA \rrlowertwocell_{G(gf)}{<-2>~~~\alpha(gf)} && GC }$
\newline\noindent for all composable 1-cells $f:A\to B$ and $g:B\to C$ in \A
\item
$\xymatrix @R1pc @C3pc  {%
      FA \ruppertwocell<5>^{1_{FA}}{<3>} \ar[dd]_{\alpha A} &
      FA \ar[dd]^{\alpha A} &&
      FA \ar[dd]_{\alpha A} \rtwocell^{1_{FA}}_{F1_A}{~~\phi^0_A} & 
      FA \ar[dd]^{\alpha A} \\
      &&=&& \\
GA \rtwocell^{1_{GA}}_{G1_A}{~~\psi^0_A} & GA && 
GA \rlowertwocell<-5>_{Gg}{<-3>~~\alpha 1_A} & GA }$ 
\newline\noindent for all $A\in\A$, where the unnamed 2-cell is the 
composite of the canonical
isomorphisms $(\alpha A)1_{FA}\cong \alpha A\cong 1_{GA}(\alpha A)$.
\end{enumerate}

Here the $\alpha A$ are called the {\em components} of the transformation, 
and the $\alpha f$ are called the {\em oplax naturality constraints}.%
\footnote{There are good reasons for regarding both the $\alpha A$ and the
$\alpha f$ as components of different type, but this point of view will 
not be used here.} 
The oplax natural transformation $\alpha$ is {\em strict} if the $\alpha f$
 are all identity 2-cells.
The case where all the $\alpha f$ are invertible is important too: then 
$\alpha$ is said to be a {\em pseudonatural} transformation. But our
main interest will be in the case where the {\em components} are identity
1-cells.

Now the bicategories, lax functors, and oplax natural transformations
do {\em not} form a 2-category. Restricting from bicategories to 2-categories
and from lax functors to 2-functors would solve two of the problems listed
below, but not the third.

\subsection*{First problem} This involves the
vertical composition of 2-cells. Given oplax natural transformations
$\alpha:F\to G$ and $\beta:G\to H$, the composite $\beta\alpha:F\to H$
has component at $A$ given by the composite $\beta A.\alpha A:FA\to HA$. 
Since composition in \B is not required to be strictly associative, there
is no reason why composition of oplax natural transformations should
be strictly associative. Notice that it's no good saying that we'll just
try to get a bicategory of bicategories instead: composition of 1-cells
{\em is} associative, it is the composition of {\em 2-cells} which is the 
problem, and in a 2-category or bicategory there is no room for weakening that;
that would require moving to a 3-dimensional structure such as a
tricategory. 

However this first problem could be solved by requiring
the objects of our desired 2-dimensional category to be 2-categories rather
than general bicategories.
If \A and \B are 2-categories, then there {\em is} a category 
${\Oplax}_l(\A,\B)$ of lax functors 
from \A to \B and oplax natural transformations, so we might hope that
it is the hom-category for a putative 2-category or bicategory of bicategories.

\subsection*{Second problem}
The second problem arises when we try to define 
composition functors 
${\Oplax}_l(\B,\C)\t{\Oplax}_l(\A,\B)\to{\Oplax}_l(\A,\C)$. 
On objects we can just use 
composition of lax functors, with which there is no problem. Given 
$\alpha:F\to G:\A\to\B$ and
$H:\B\to\C$ we need to be able to define $H\alpha:HF\to HG$. This would involve
components $H\alpha A:HFA\to HGA$ and the obvious choice is simply to
apply $H$ to $\alpha A:FA\to GA$. We now need oplax naturality maps
$$\xymatrix{
HFA \ar[rr]^{H\alpha A} \ar[d]_{HFf} &{}& HGA \ar[d]^{HGf} \\
HFB \ar[rr]_{H\alpha B} & {}\utwocell\omit{(H\alpha) f} & HGB }$$
but here we cannot simply apply $H$ to $\alpha f:\alpha B.Ff\to Gf.\alpha A$,
since this would give a map $H(\alpha B.Ff)\to H(Gf.\alpha A)$. We could
use the lax functoriality constraint of H to get a map 
$H\alpha B.HFf\to H(Gf.\alpha A)$, but the codomain would still be wrong.
We could fix this in the same way if $H$ were not just a lax functor but a 
homomorphism, but problems would still remain.

\subsection*{Third problem}

Even if we restricted to 2-categories and 2-functors there would still
be a problem,
once again involving the definition of composition. Given 
$\alpha:F\to G:\A\to\B$ and $\beta:H\to K:\B\to\C$ 
we can define $H\alpha$ and $\beta F$ (and $K\alpha$, and $\beta G$)
and so we do obtain a sesquicategory \Oplax. But the middle four
interchange law, asserting the equality of the vertical composites
$$\xymatrix @R1pc  {
\A \rtwocell^{F}_{G}{\alpha} & \B \ar@/^1pc/[r]^{H} & \C &&
\A \ar@/^1pc/[r]^{F} & \B \rtwocell^{H}_{K}{\beta} & \C \\
\A \ar@/_1pc/[r]_{G}  & \B \rtwocell^{H}_{K}{\beta} & \C &=&
\A \rtwocell^{F}_{G}{\alpha} & \B \ar@/_1pc/[r]_{K} & \C }$$
need not hold. It does hold when $\beta$ is 2-natural (strictly natural)
but not in general; and there is once
again no room to weaken this law: that would require a third dimension.

\subsection*{Solution}
There is, however, another way of restricting the 2-cells which 
solves all three problems. We still allow arbitrary bicategories as objects,
and arbitrary lax functors as 1-cells, but an oplax natural transformation 
$\alpha:F\to G$ is allowed as a 2-cell only if the components $\alpha A$
are all identity 1-cells. Such an $\alpha$ will be called an {\em icon}. 

If the components of $\alpha$ are to be identities, then $F$ and $G$ 
have to agree on objects. In fact we then omit the components
altogether, so as to avoid any problems with the identities in \B being
non-strict. We shall now formulate the precise definition. We use
$M$ to denote a composition map such as $\A(B,C)\t\A(A,B)\to\A(A,C)$,
except that to save space we shall write this as $M_{A,B,C}:\A_2\to\A_1$.
Similarly we write $j_A:1\to\A_1$ for the map $1\to\A(A,A)$ picking out
the identity on $A$. Recall also that we write $\phi$ and $\phi^0$ for the 
lax-functoriality
constraints of a lax functor $F$, and similarly $\psi$ and $\psi^0$ for $G$.

\begin{definition}
If $F,G:\A\to\B$ are lax functors of bicategories, and $FA=GA$ for all
objects $A\in\A$, an {\em icon} from $F$ to $G$ consists of a natural
transformation
$$\xymatrix @C3pc {
{}\llap{\A(A,B}) \rtwocell^{F}_{G}{\alpha} & {}\B\rlap{(FA,FB)} }$$
for all objects $A,B\in\A$, subject to the following two conditions:
\begin{enumerate}[({I}CON1)]
\item 
$$\xymatrix @R1pc @C4pc {%
\A_2 \rtwocell^{F\t F}_{G\t G}{~~~\alpha\t\alpha} \ar[dd]_(0.6){M_{A,B,C}} &
\B_2 \ar[dd]^(0.6){M_{FA,FB,FC}} &&
\A_2 \ruppertwocell<2>^{F\t F}{<4>\phi} \ar[dd]_(0.6){M_{A,B,C}} & 
\B_2 \ar[dd]^(0.6){M_{FA,FB,FC}} \\
&&=&& \\
\A_1 \rlowertwocell<-2>_{G}{<-4>\psi} & \B_1 &&
\A_1 \rtwocell^{F}_{G}{\alpha} & \B_1
}$$
\item
$$\xymatrix @R2pc {%
& 1 \ar[ddr]^{j_{FA}} \ar[ddl]_{j_A} \ddtwocell\omit{\phi^0} &&&
& 1 \ar[ddr]^{j_{FA}} \ar[ddl]_{j_A} \ddtwocell\omit{\psi^0} \\ 
&&&=&&& \\
\A_1 \rrtwocell^{F}_{G}{\alpha} && \B_1 &&
\A_1 \ar@/_0.5pc/[rr]_{G} && \B_1 }$$
\end{enumerate}
\end{definition}

Here (ON0) corresponds to the naturality of $\alpha:F\to G:\A(A,B)\to\B(FA,FB)$,
while (ON1) and (ON2) correspond to (ICON1) and (ICON2), respectively.

There is now no problem in obtaining a 2-category \bicat of bicategories,
lax functors, and icons. For bicategories \A and \B, the hom-category
$\bicat(\A,\B)$ has lax functors from \A to \B as objects, and icons
as morphisms. Given icons $\alpha:F\to G$ and $\beta:G\to H$, we have
$FA=GA=HA$ for all $A\in\A$, and the composite $\beta\alpha:F\to H$ is
defined using the composite natural transformations
$$\xymatrix{%
\A(A,B) \ruppertwocell^{F}{\alpha} \ar[r]^(0.4){G} \rlowertwocell_{H}{\beta} &
\B(FA,FB) }$$
and composition of natural transformations is of course associative,
so that we do indeed have a category $\bicat(\A,\B)$.

Given an icon $\alpha:F\to G:\A\to\B$ and lax functors $H:\B\to\C$
and $K:\D\to\A$, the composite $H\alpha K:HFK\to HGK$ is defined using
the natural transformations
$$\xymatrix {%
\D(D,E) \ar[r]^-{K} & \A(KD,KE) \rrtwocell^{F}_{G}{\alpha} && 
\B(FKD,FKE) \ar[r]^-{H} & \C(HFKD,HFKE) }$$
and the middle-four interchange law holds, so that we do indeed have
a 2-category \bicat:

\begin{theorem}
There is a 2-category \bicat of bicategories, lax functors, and icons,
with the usual composition laws.
\end{theorem}

\section{Bicategories as many-object monoidal categories}

Every monoidal category \V determines a one-object bicategory $\Sigma\V$.
If we write $*$ for the unique object of $\Sigma\V$, then the hom-category
$\Sigma\V(*,*)$ is the category \V, composition is given by 
$\ot:\V\t\V\to\V$, and the identity on $*$ by the unit object $I$. A
 monoidal functor $\V\to\W$ is the same thing as a lax functor
$\Sigma\V\to\Sigma\W$. But now if $F,G:\V\to\W$ are
 monoidal functors, and $\Sigma F,\Sigma G:\Sigma\V\to\Sigma\W$ 
the corresponding lax functors, what would it be to give an oplax
natural transformation from $\Sigma F$ to $\Sigma G$? Since 
$\Sigma\V$ has one object, this would have a single component, given
by an object $W\in\W$. Oplax naturality would mean giving, for every
$V\in\V$, a morphism $W\ot FV\to GV\ot W$ satisfying various conditions.
This is rather more general than a monoidal transformation from
$F$ to $G$ (this fact was particularly significant in 
\cite{Cheng-Gurski}): a monoidal transformation is a natural transformation 
$\alpha:F\to G$ which is suitably compatible with the  monoidal
structures of $F$ and $G$, and there is certainly no role for the 
object $W$ appearing in an oplax natural transformation $\Sigma F\to\Sigma G$.
But an {\em icon} from $\Sigma F$ to $\Sigma G$ is exactly a monoidal
natural transformation from $F$ to $G$.

We conclude:

\begin{theorem}
The 2-category \MonCat of monoidal categories, monoidal functors, 
and monoidal natural transformations can be identified with the full
sub-2-category of \bicat consisting of the one-object bicategories.
\end{theorem}

There is also a corresponding result where we use homomorphisms in 
place of lax functors and strong monoidal functors in place of monoidal
functors.

In some contexts the objects of a bicategory are the most important level
of structure. But in other contexts they are really just a way of parametrizing
the various 1-cells and 2-cells. The  bicategory \Mod
 of rings, (bi)modules,
and module homomorphisms is a good example of the latter situation. A morphism
from $R$ to $S$ is a left $R$-, right $S$-bimodule, and a 2-cell is a 
homomorphism of such bimodules. Much of the time the rings are of rather
secondary importance. In situations like this, one can think of a 
bicategory as a ``many-object'' monoidal category. Our 2-category \bicat
is particularly well-suited to dealing with bicategories thought of in
this light. In the remainder of this section we look at equivalences 
in \bicat as an example of the sort of thing we have in mind.

In any 2-category one has the notion of equivalence: a 1-cell $f:A\to B$
with a $g:B\to A$ satisfying $gf\cong 1$ and $fg\cong 1$. In particular,
one could apply this to the case of \bicat. First we characterize the
invertible 2-cells:

\begin{proposition}
An icon $\alpha$ from $F$ to $G$ is invertible if and only if for all
objects $A$ and $B$ the natural
transformation
$$\xymatrix{
{}\llap{\A(A,B}) \rtwocell^{F}_{G}{\alpha} & {}\B\rlap{(FA,FB)} }$$
is invertible. This in turn is saying that the oplax 
natural transformation is pseudonatural.
\end{proposition}

\begin{proposition}
A lax functor $F:\A\to\B$ is an equivalence in \bicat if and only if
(i) $F$ is bijective on objects, (ii) each $F_{A,B}:\A(A,B)\to\B(FA,FB)$
is an equivalence of categories, and (iii) $F$ is a homomorphism.
\end{proposition}

\proof
 An equivalence in \bicat 
consists of a lax functor $F:\A\to\B$ for which there is another lax
functor $G:\B\to\A$ with invertible icons $GF\cong 1$ and $FG\cong 1$. 
The existence of these icons already forces the object-parts of $F$ and
$G$ to be mutually inverse. In order to analyze the remaining structure
of an equivalence in \bicat, we may as well suppose that \A and \B have 
the same
objects, and that $F$ and $G$ act as the identity on objects. For each 
pair $A,B\in\A$ of objects, we then have functors $F_{A,B}:\A(A,B)\to\B(A,B)$
and $G_{A,B}:\B(A,B)\to\A(A,B)$, and the invertible icons provide natural
isomorphisms $G_{A,B}F_{A,B}\cong 1$ and $F_{A,B}G_{A,B}\cong 1$, thus 
each $F_{A,B}$ is indeed an equivalence. 

We still need to show that $F$ is a homomorphism. Consider the invertible
icon $\beta:GF\cong 1$. For 1-cells $f:A\to B$ and $g:B\to C$ in \A, we
have the comparisons $\phi_{g,f}:Fg.Ff\to F(gf)$ and 
$\psi_{Fg,Ff}:GFg.GFf\to G(Fg.Ff)$. Part of the condition of $\beta$ to 
be an icon is the commutativity of 
$$\xymatrix{
GFg.GFf \ar[r]^{\psi_{Fg,Ff}} \ar[d]_{\beta.\beta} & 
G(Fg.Ff) \ar[r]^{G\phi_{g,f}} & GF(gf) \ar[d]^{\beta} \\
gf \ar@{=}[rr] && gf, }$$
in which the vertical maps are invertible, and so the composite across
the top is so too. Thus $\psi_{Fg,Ff}$ is split monic for all $f$ and $g$;
since $F$ is an equivalence on hom-categories, it follows that 
$\psi_{h,k}:Gh.Gk\to G(hk)$ is split monic for all $h:A\to B$ and $k:B\to C$
in \B. By symmetry, each $\phi_{g,f}$ is also split monic. On the other 
hand, $G\phi_{g,f}$ is split epi, and so since $G$ is an equivalence
on hom-categories, each $\phi_{g,f}$ is split epi, and so invertible. 
The argument that each $\phi^0_0:1_A\to F1_A$ is invertible is similar
but easier, and so we conclude that $F$ is indeed a homomorphism.

Now we turn to the converse: any lax functor satisfying the three conditions
is an equivalence. It is
well-known that for any homomorphism of bicategories $F:\A\to\B$ which is 
is an equivalence on hom-categories and surjective on objects up to 
equivalence, there is a homomorphism $G:\B\to\A$ with the composites 
$GF$ and $FG$ each pseudonaturally equivalent to the relevant identity
homomorphism; such an $F$ is called a biequivalence. More precisely, given 
such an $F:\A\to\B$ and a choice,
for each $B\in\B$ of an object $GB\in\A$ and an equivalence 
$eB:FGB\to B$, one can extend $G$ to a homomorphism of bicategories,
so that the $eB$ become the components of a pseudonatural 
transformation. Furthermore, since $F$ is an equivalence on hom-categories,
we can find for each $A\in\A$ a morphism  $dA:GFA\to A$ in \A with 
$FdA:FGFA\to FA$ isomorphic to $eFA:FGFA\to FA$, and any 
such choice can be extended to a pseudonatural equivalence $d:GF\to 1$. 

In our case, $F$ is bijective on objects, so we can
choose $GB$ so that $FGB=B$. Thus we obtain a homomorphism $G:\B\to\A$
and an invertible icon $FG\cong 1$. Now the $eB$ of the previous paragraph
is the identity, so we can take $dA:GFA\to A$ also to be the identity, and
then the resulting $d:GF\to 1$ becomes our invertible icon $GF\cong 1$.
\endproof

It is really the biequivalences that provide the general notion of 
``sameness'' for bicategories; for example Gabriel-Ulmer duality 
asserts the biequivalence, in this sense, of the 2-category \Lex
of categories with finite limits, and a certain 2-category \lfp of locally
finitely presentable categories. This biequivalence involves genuine 
content at the object-level, and is certainly not an equivalence 
in \bicat. 

On the other hand, there are some important examples, where ``nothing
happens'' with the objects. A good example is the theorem that every
bicategory is biequivalent to a 2-category (``every bicategory can be
made strict'') \cite{MacLane-Pare}. Here the lack of strictness has nothing 
to do with the
objects, and they can be left unchanged; the problem is rather with the
1-cells. Thus it {\em is} the case that every bicategory is equivalent in
\bicat to a strict one, and we can conclude (see also \cite{2-nerve}):

\begin{theorem}
The 2-category \bicat is biequivalent to the full sub-2-category consisting
of the strict bicategories (the 2-categories).
\end{theorem}

\section{Costrictness}

In this section we describe an abstract characterization of the icons 
among all oplax natural transformations. We restrict, for convenience,
to the case of 2-categories and 2-functors, although a similar analysis
could be given involving bicategories and homomorphisms (not lax functors).
The formal notion of strictness in a sesquicategory was defined in 
Section~\ref{sect:background}; here we work through the idea more gently
in the specific case of the sesquicategory \Oplax of 2-categories,
2-functors, and oplax natural transformations.

Consider then a diagram of 2-categories, 2-functors, and oplax natural
transformations
$$\xymatrix{
\A \rtwocell^{F}_{G}{\alpha} & \B \rtwocell^{H}_{K}{\beta} & \C. }$$
As described in Section~\ref{sect:background} above, we can form 
$\beta G:HG\to KG$ and $H\alpha:HF\to HG$, and their composite
$\beta G.H\alpha:HF\to KG$; and similarly the composite
$K\alpha.\beta F:HF\to KG$. The middle-four interchange law would
require these composites to be the same, but in general they are not. 

We can nonetheless, consider situations under which they are the same. 
We start with the following question:

\begin{quotation}
\noindent
  For which $\beta:H\to K:\B\to\C$ is it the case that middle-four
  interchange holds for {\em all} $\alpha:F\to G:\A\to\B$?
\end{quotation}

In particular, it would have to hold when \A is the terminal 2-category 1. 
Then $F$ and $G$ are just objects of \B, and $\alpha:F\to G$ a morphism,
and middle-four interchange says precisely that the components 
$\beta B:HB\to KB$ are strictly natural. Notice that this does not yet
imply that $\beta$ is itself strictly natural, for we might have strictly
natural components yet still choose non-trivial oplax naturality maps.

But now let \A be the arrow 2-category $\two$. A 2-functor $\two\to\B$
is just a morphism $f:A\to B$ in \B. An oplax natural transformation
between two such 2-functors $\two\to\B$, say $f:A\to B$ and $g:C\to D$,
is a square 
$$\xymatrix @R1pc {%
A \ar[r]^{f} \ar[dd]_{u} & B \ar[dd]^{v} \\ 
{} \rtwocell\omit{\phi} & {} \\
C \ar[r]_{g} & D. }$$
Now let $G:\two\to\B$ 
correspond to a map $f:A\to B$ in \B, let $F:\two\to\B$ correspond to 
$1:A\to A$, and let $\alpha:F\to G$ correspond to the strict oplax natural 
transformation with components $1:A\to A$ and $f:A\to B$, as in 
$$\xymatrix {%
A \ar@{=}[r] \ar@{=}[d] & A \ar[d]^{f} \\
A \ar[r]_{f} & B.}$$
Then the 
equation $\beta G.H\alpha=K\alpha.\beta F$ becomes
$$\xymatrix @R1pc {
HA \ar[r]^{1} \ar[dd]_{1} & HA \ar[dd]^{Hf} &&
HA \ar[r]^1 \ar@{=}[dd] & HA \ar@{=}[dd] \\
&&& {}\rtwocell\omit{~~\beta1_A} & {}\\
HA \ar[r]^{Hf} \ar@{=}[dd] &
HB \ar@{=}[dd] &=&
KA \ar[r]^1 \ar[dd]_{1} & KA \ar[dd]^{Kf}  \\
 {}\rtwocell\omit{~\beta f} & {} \\
KA \ar[r]_{Kf} & KB && KA \ar[r]_{Kf} & KB }$$
but $\beta1_A$ is the identity, and so this says that $\beta f$ is
the identity. Since $f$ was arbitrary,
this says that $\beta$ really is strict. 

Conversely, if $\beta$ is strict, then clearly the middle-four interchange
law will hold for all $\alpha$. But now we have an internal (representable)
notion of strictness, and so we can consider the dual:

\begin{quotation}
\noindent  For which $\alpha:F\to G:\A\to\B$ is it the case that middle-four
  interchange holds for {\em all} $\beta:H\to K:\B\to\C$?
\end{quotation}

Formally, one could define strictness in any sesquicategory. Then our
original notion would be strictness in $\Oplax$, while the question
just posed asks which $\alpha$ are strict in $\Oplax\op$. 
Such an $\alpha$ will be called {\em costrict}. We shall show that 
$\alpha$ is a costrict if and only if it is an icon.

Let $\alpha$ and
$\beta$ be as above, with $\alpha$ an icon. We must show that 
$\beta G.H\alpha=K\alpha.\beta F$, but this amounts to showing that
these composites have the same component for each $A\in\A$, and the
same oplax naturality map for each $f:A\to B$ in \A. Now $FA=GA$ and
$\alpha A$ is the identity, so 
$(\beta G.H\alpha)A=\beta GA.H\alpha A=\beta GA=\beta FA=K\alpha A.\beta FA
=(K\alpha.\beta F)A$ and the components agree. As for the second part,
it asserts the equality of pasting composites
$$\xymatrix{%
HFA \ar[r]^{HFf} \ar@{=}[dd] & HFB \ar@{=}[dd] &&
HFA \ar[r]^{HFf} \ar[dd]_{\beta FA} & HFB \ar[dd]^{\beta FB} \\
{} \rtwocell\omit{~~H\alpha f} & {} && {}\rtwocell\omit{~~\beta Ff} & {} \\
HGA \ar[r]^{HGf} \ar[dd]_{\beta GA} & HGB \ar[dd]^{\beta GB} &=&
KFA \ar[r]^{KFf} \ar@{=}[dd] & KFB \ar@{=}[dd] \\
{} \rtwocell\omit{~~\beta Gf} & {} && {}\rtwocell\omit{~~K\alpha f} & {} \\
KGA \ar[r]_{KGf} & KGB && KGA \ar[r]_{KGf} & KGB }$$
but this follows by (ON0) for the oplax natural transformation $\beta$.

Suppose conversely, that $\alpha$ is strict.
We shall show that it must be an icon. We do this by constructing a 
2-category \C with 2-functors $H,K:\B\to\C$ and an oplax natural transformation
$\beta:H\to K$ with the property that if $g:B\to C$ is a morphism in \B
with respect to which $\beta$ is strictly natural, then $g$ must be an
identity. Since $\beta$ is by assumption strictly natural with respect
to each component $\alpha A:FA\to GA$, it will follow that these components
are identities, and so that $\alpha$ is an icon.

How do we construct such a \C? Among all the 2-categories \C equipped
with 2-functors $H,K:\B\to\C$ and an oplax natural transformation 
$\beta:H\to K$, there is a universal one, given by the (lax) Gray
tensor product $\two\ot\B$ of the arrow category $\two=\{0<1\}$ and \B.
The objects of $\two\ot\B$
are pairs $(i,B)$ where $i\in\two$ and $B\in\B$. The morphisms
of $\two\ot\B$ are freely generated by $(!,B):(0,B)\to(1,B)$ for each $B\in\B$,
and $(i,g):(i,B)\to(i,C)$ for each $i\in\two$ and $g:B\to C$ in \B;
subject to the relations that $(i,hg)=(i,h)(i,g)$ and $(i,1_B)=1_{(i,B)}$.
The 2-cells are generated by $(i,\beta):(i,g)\to(i,h)$ for all
2-cells $\beta:g\to h$ in \B, and 
$$\xymatrix @R1pc {%
(0,B) \ar[r]^{(0,g)} \ar[dd]_{(!,B)} & (0,C) \ar[dd]^{(!,C)} \\
{} \rtwocell\omit{~~(!,g)} & {} \\
(1,B) \ar[r]_{(1,g)} & (1,C) }$$
for all 1-cells $g:B\to C$, subject to certain relations. Which ones?
The 2-functors $H$ and $K$ are defined so as to send $B\in\B$ to 
$(0,B)$ and to $(1,B)$, respectively, and similarly on morphisms and 2-cells.
We want the $(!,B):(0,B)\to(0,B)$ to be the components $HB\to KB$ of
an oplax natural transformation, with oplax naturality constraints 
$(!,g)$. The relations are exactly the ones which make this work. 
(See \cite{Gray} or \cite{Street:categorical-structures} for details.) The
important thing for us is that $(!,g)$ is not an identity 2-cell unless 
$g$ is an identity 1-cell, but in fact we can see this without knowing
exactly what the 2-cells are, because the domain and codomain 1-cells
$(!,C)(0,g)$ and $(1,g)(!,B)$ of $(!,g)$ are different unless $g$ is
an identity 1-cell. 

This now proves:

\begin{theorem}
An oplax natural transformation between 2-functors is costrict if 
and only if it is an icon.
\end{theorem}

\begin{remark}
We have defined the notion of strictness in any sesquicategory. 
Just as one-object 2-categories can be identified with strict 
monoidal categories, one-object sesquicategories can be identified
with structures called {\em premonoidal categories}
\cite{Power-Robinson:premonoidal}. A morphism in a premonoidal 
category is called {\em central} \cite{Power-Robinson:premonoidal} 
if the corresponding 2-cell is both strict and costrict. In
the sesquicategory of 2-categories, 2-functors, and oplax natural
transformations, the only 2-cells which are strict and costrict are
the identities, but in general it is possible to have non-identity
central morphisms; see \cite{Power-Robinson:premonoidal}.
\end{remark}

\section{Alternative approaches}

In this section we give two alternative approaches to icons, one using
double categories, the other using 2-dimensional universal algebra.

\subsection{Bicategories as pseudo double categories}

A double category involves two categories with the same set of objects,
and with the morphisms called vertical and horizontal, respectively. There are
also ``squares'' which have a vertical domain, a vertical codomain, a 
horizontal domain, and a horizontal codomain; and which can be composed
either horizontally or vertically. For instance in 
$$\xymatrix{
A \ar[r]^{h} \ar[d]_{v} \ar@{}[dr]|{\alpha} & B \ar[d]^{v'} \\
C \ar[r]_{h'} & D }$$
$h:A\to B$ and $h':C\to D$ are horizontal arrows, and are the horizontal
domain and horizontal codomain of the square $\alpha$, while 
$v:A\to C$ and $v':B\to D$ are vertical arrows, and are the vertical 
domain and vertical codomain of the square.

A 2-category can be thought of as a
double category with no non-identity horizontal arrows, so that all ``squares''
collapse to the usual shape of a 2-cell in a 2-category, as in 
$$\xymatrix{
A \dtwocell^{h'}_{h}{^\alpha} \\ B }$$
A pseudo double
category is a slight generalization of a double category, in which the 
vertical, but not the horizontal, structure is weakened. A pseudo
double category with no non-identity horizontal cells is the same thing as a 
bicategory. 

In \cite{GrandisPare05}, Grandis and Par\'e defined a 2-category 
\LaxDbl whose objects were pseudo double categories, whose morphisms
are lax functors, and ``horizontal transformations'' as 2-cells. 
Here a horizontal transformation involves components which are horizontal
morphisms; these are only weakly natural. If we restrict the objects
of \LaxDbl to the bicategories (seen as pseudo double categories with
no non-identity horizontal cells), then such horizontal transformations
are forced to have identity components, and in fact turn out to be icons.

\subsection{Bicategories as algebras for a 2-monad}

In \cite{2-nerve} an alternative approach to icons was also given.
There a 2-category \CG of \Cat-enriched graphs was defined. A \Cat-graph
\G has vertices $A,B,C,\ldots$, and ``hom-categories'' $\G(A,B)$ for all 
vertices $A$ and $B$. A homomorphism $F:\G\to\HH$ of \Cat-graphs consists 
of a function $A\mapsto FA$ sending vertices of \G to vertices of \HH, and
a functor $F:\G(A,B)\to\HH(FA,FB)$ for all vertices $A$ and $B$. The 
\Cat-graphs and their homomorphisms form a category \cite{Wolff}, but
it is also possible to make it into a 2-category. A 2-cell $F\to F'$
is possible only if $F$ and $F'$ agree on vertices; it then consists of
a natural transformation 
$$\xymatrix @C3pc {
{}\llap{\G(A,B}) \rtwocell^{F}_{F'}{\alpha} & \HH\rlap{(FA,FB)} }$$
for all vertices $A$ and $B$, with no further conditions. 

This 2-category is locally finitely presentable in the sense of 
\cite{Kelly-amiens}, and one can give a presentation, in the sense of 
\cite{KP}, of a finitary 2-monad on \CG whose algebras are the bicategories,
whose strict morphisms are the strict homomorphisms, and whose algebra
2-cells are the icons. The pseudomorphisms of algebras are the homomorphisms,
and so the 2-category \homm is the 2-category \talg of strict algebras, 
pseudomorphisms, and algebra 2-cells, which is the main object of study
in \cite{BKP}. It follows that \homm has
bicategorical limits and colimits, and strict products, inserters, 
and equifiers. All this was done in  \cite{2-nerve}.

It is not hard to verify similarly that the lax morphisms of algebras are
the lax functors. Thus
the 2-category \bicat of bicategories, lax functors, and icons is the 2-category
\talgl of strict algebras, lax morphisms, and algebra 2-cells, which is
the main object of study of \cite{talgl}, and so by \cite{talgl} 
\bicat has oplax limits, and in particular products, cotensor products,
Eilenberg-Moore objects of comonads, {\em some} inserters and equifiers,
and any limit of strict morphisms.

\section{Applications}

\subsection{2-nerves}

There is a functor from the category \cat of categories to the category
$[\Delta\op,\Set]$ of simplicial sets, which sends a category to its 
{\em nerve}. Concretely,
the 0-simplices of the nerve are the objects, the 1-simplices the morphisms,
the 2-simplices the composable pairs, and so on. But the slick way to define
this is to regard $\Delta$ as the full subcategory of \cat consisting of 
the finite ordinals (the finite totally ordered sets, regarded
as categories), and then use the 
inclusion $J:\Delta\to\cat$ to define the nerve $NC$ of a category $C$
as the simplicial set $\cat(J-,C)$ sending an object $\n\in\Delta$ to 
the set $\cat(J\n,C)$ of functors from $J\n$ to $C$.

More generally, if \C is a category enriched in \V, and $J:\A\to\C$ is 
a \V-functor, each object $C$ determines a \V-functor $\C(J-,C)$ sending
$A\in\A$ to the hom-object $\C(JA,C)\in\V$; now there is a \V-functor
$\C\to[\A\op,\V]$ sending $C$ to $\C(J-,C)$.

We shall apply this in the special case where $\V=\Cat$, so a \V-category
is a 2-category. We take \C to be the 2-category \nhom of bicategories,
normal homomorphisms, and icons, and we take \A to be $\Delta$, seen 
as a 2-category with no non-identity 2-cells. We have the inclusion
of $\Delta$ in $\cat$, but then any category can be regarded as a 
bicategory with no non-identity 2-cells, and so we get an inclusion
$J:\Delta\to\nhom$, which turns out to be fully faithful. The induced
\V-functor $\nhom\to[\Delta\op,\Cat]$ sends a bicategory $B$ to a
simplicial object $\nhom(J-,B)$ in \Cat, called the {\em 2-nerve} of 
the bicategory. See \cite{2-nerve} for details.

\subsection{Lawvere theories}

A Lawvere theory is a category \T with finite products, equipped with 
a functor $\sS\op\to\T$ from the opposite of the category of finite 
sets, which is bijective on objects and preserves finite products.
In other words, \T is a category with finite products whose objects have 
the form $X^n$ for $n\in\N$, and with $X^m\t X^n=X^{m+n}$. A morphism of
Lawvere theories from $E:\sS\op\to\T$ to $E':\sS\op\to\T'$ is a functor 
$M:\T\to\T'$ with $ME=E'$; then $M$ necessarily preserves finite products. 

There is a 2-dimensional version of Lawvere theory (which is a special case 
of the more general
enriched categorical version \cite{Power:enriched-lawvere}). In place of \sS 
we have the 2-category 
\cC of finitely presentable categories, and now \T is required to have
{\em cotensors} by objects of \cC (finite cotensors for short) and a 
bijective-on-objects finite-cotensor-preserving 2-functor $E:\cC\op\to\T$
Once again, a morphism of theories from $E:\cC\op\to\T$ to $E':\cC\op\to\T'$
is a 2-functor $M:\T\to\T'$ satisfying $ME=E'$, and once again such an $M$
necessarily preserves finite cotensors. But there is now an ingredient 
not present in the general enriched version: the category of \V-enriched
Lawvere theories is just a category, there is no enrichment over \V; but
the category of 2-dimensional Lawvere theories does enrich over \Cat to 
give a 2-category. This corresponds precisely to the fact that there is 
a 2-category of finitary 2-monads on \Cat (studied in \cite{KP} and 
\cite{mnd}). Given morphisms $M,M':\T\to\T'$ as above, a 2-cell
from $M$ to $M'$ is an oplax natural transformation $\phi:M\to M'$ whose
restriction along $E$ is the identity transformation on $E'$. This 
condition implies in particular that $\phi$ is an icon, and this means
that we do indeed get a 2-category of 2-dimensional Lawvere theories. 
This is important when one wishes to consider {\em flexibility} of theories.

Furthermore, the morphisms $M$ and $M'$ induce 2-functors 
$M,M':\Mod(\T')\to\Mod(\T)$ between the 2-categories of models, and 
now composition with $\phi$ induces a 2-natural transformation $M\to M'$.

See \cite{2dLawvere} for more details.

\subsection{Bundles}

We have already seen that the 2-category structure of \bicat can be 
used to obtain a notion of equivalence of bicategories. It can also
be used to obtain a notion of fibration. 

Briefly, if $p:A\to B$ is a morphism in a 2-category \K, we can define
what it means for $p$ to be a fibration. First of all, for a morphism
$a:X\to A$, a 2-cell $\alpha:a'\to a$ is said to be {\em $p$-cartesian} 
if for each $c:X\to A$ and each pair $(\gamma:pc\to pa',\delta:c\to a)$ 
with $p\alpha.\gamma=p\delta$ there is a unique
$\bar{\gamma}:c\to a'$ with $p\bar{\gamma}=\gamma$ and $\alpha\gamma=\delta$,
as in the diagram below.
$$\xymatrix {%
c \ar[drr]^{\delta} \ar@{.>}[dr]_{\bar{\gamma}} \\
& a' \ar[r]_{\alpha} & a \\
pc \ar[drr]^{p\delta} \ar[dr]_{\gamma} \\
& pa' \ar[r]_{p\alpha} & pa }$$
Then $p$ is said to be a {\em fibration} if (i) for each $a:X\to A$, $b:X\to B$,
and $\beta:b\to pa$ there exists a $p$-cartesian $\alpha:a'\to a$ with 
$pa'=b$ and $p\alpha=\beta$, and (ii) if $\alpha:a'\to a$ is $p$-cartesian,
then $\alpha x:a'x\to ax$ is $p$-cartesian for any $x:Y\to X$. 

One can now consider fibrations in the 2-category \bicat of bicategories,
normal homomorphisms, and icons. A fibration will be a homomorphism 
$P:\A\to\B$ of bicategories, with for each $f:A\to A'$ in \A and each
2-cell $\beta:g\to Pf$ in \B a 2-cell $\bar{\beta}:\bar{g}\to f$ in \A
with $P\bar{\beta}=\beta$, as in the diagram
$$\xymatrix @R3pc {
A \ar[r]^{f} \ar@/_1pc/@{.>}[r]_{\bar{g}} \rtwocell\omit{^<1>\bar{\beta}} & A' \\
PA \ar[r]^{Pf} \ar@/_1pc/[r]_{g} \rtwocell\omit{^<1>\beta} & PA' }$$
and with the liftings $\bar{\beta}$ required to satisfy various properties. 

There are various special cases. In particular, one could consider 
fibration not just in \bicat, but in the 2-category \shom of bicategories,
{\em strict} homomorphisms, and icons. This imposes further compatibility
conditions on the liftings $\bar{\beta}$. A still greater restriction 
would be to ask not just for a fibration in \shom but a split fibration. 
This means that one makes chosen liftings, which are themselves functorial.
This situation is being studied in \cite{HL} under the name {\em bundles of
bicategories} after it arose in \cite{Hess-Parent-Scott}; see \cite{HL} 
for examples and further discussion.


\bibliographystyle{plain}

\end{document}